\documentclass[preprint,12pt]{elsarticle}
  \usepackage{amsmath} 
  \usepackage{amssymb} 
  \usepackage{MnSymbol}
  \usepackage{extarrows}
  \usepackage{enumerate} 
  
  \usepackage{bm} 
  
  \usepackage{graphicx} 
  \usepackage{subfigure} 
  \usepackage{multirow} 
  
  \usepackage[colorlinks = true, citecolor = blue,  urlcolor = red]{hyperref}
  \newtheorem{theorem}{Theorem}[section]
  \newtheorem{definition}{Definition}[section]
  
  \newtheorem{lemma}{Lemma}[section]

  \newtheorem{corollary}{Corollary}[section]
  \newtheorem{remark}{Remark}[section]
  \renewcommand{\P}{\mathbb{P}}
  \newcommand{\E}{\mathbb{E}}
  \newcommand{\V}{\text{Var}}
  \newcommand{\pf}{\textbf{Proof: }}
  \newcommand{\e}{\hfill$\blacksquare$}
  
  \newcommand{\R}{\mathbb{R}}
  \newcommand{\KL}{\text{KL}}

  \allowdisplaybreaks[4] 
  \numberwithin{equation}{section} 

\begin{document}

 \begin{frontmatter}
\title{New Upper bounds for KL-divergence Based on Integral Norms}

	\author[label1]{Liu-Quan Yao}
\ead{yaoliuquan20@mails.ucas.ac.cn}
\author[label2]{Song-Hao Liu}
\ead{liusonghao_cuhk@link.cuhk.edu.hk}

\address[label1]{Academy of Mathematics and Systems Science, Chinese Academy of Sciences,\\
	Beijing,
	100190,
	China}

\address[label2]{Department of Statistics and Data Science,
	Southern University of Science and Technology,\\
	Shenzhen,
	518055,
	Guangdong,
	China}


%

\begin{abstract}
In this paper, some new upper bounds for Kullback-Leibler divergence(KL-divergence) based on $L^1, L^2$ and $L^\infty$ norms of density functions are discussed. Our findings unveil  that the convergence in KL-divergence sense  sandwiches between the convergence of density functions in terms of $L^1$ and $L^2$ norms.   Furthermore, we endeavor to apply our newly derived upper bounds to the analysis of the rate theorem of  the entropic conditional central limit theorem.
\end{abstract}

%
%

\end{frontmatter}




	\section{Introduction}
Pinsker's inequality, a well-known inequality, shows that the KL-divergence effectively governs the total variation distance between two random variables, thus establishing a profound connection between the information-entropy theorem and classical probability theory. Specifically, considering two probability measures denoted by $P$ and $Q$, defined on a common measurable space $(\Omega, \mathcal{F})$, the Pinsker's inequality (see in \cite{Brillinger1964InformationAI} for example)  states that:
\begin{equation}\label{Pinsker's inequality}
\KL(P\| Q)\ge \frac{1}{2}\|P-Q\|_{TV}^2,
\end{equation}
where the KL-divergence is defined as $\KL(P\|Q):=\int_\Omega \ln(\frac{dP}{dQ})dP$, while the total variation distance is defined as $\|P-Q\|_{TV}:=2\sup_{A\in\mathcal{F}} |P(A)-Q(A)|$.

On the one hand, it is evident that if $P$ and $Q$ are discrete probability measures over a countably infinite set $B$, we obtain:

\begin{equation}
\|P-Q\|_{TV}=2\sup_{A\in\mathcal{F}} |P(A)-Q(A)|=\sum_{x\in B} |P(x)-Q(x)|.
\end{equation}

On the other hand, if $P$ and $Q$ are continuous probability measures on the real line $\R$, with densities $p$ and $q$, respectively, we have:

\begin{equation}
\|P-Q\|_{TV}=2\sup_{A\in\mathcal{F}} |P(A)-Q(A)|=\int_\R |p(x)-q(x)|dx=\|p-q\|_1.
\end{equation}

This relationship of control can be easily extended to the sequential case that as $n\to \infty$, 
$$\KL(P_n\|Q)\to 0\Rightarrow \|P_n-Q\|_{TV}\to 0.$$
 As a result, the use of entropy methods in classical probability theorems, such as the entropic central limit theorem, has garnered significant attention in numerous research works (\cite{Bar1986}, \cite{Carlen1991}, \cite{Johnsonbook},   \cite{Esseenbound} and \cite{ECCLT} \textit{et al.}).
	
	In recent researches,  the upper bounds for KL-divergence was discussed frequently, since KL-divergence is a common measure for random variables in practical applications (\cite{SKL-diagnostic}, \cite{SKL-Capacity1} and  \cite{SKL-Wasserstein} \textit{et al.}). One of the important series of studies for the upper bounds is the reverse Pinsker's inequality. 
 Sason (2015) \cite{reversePin2015} built a reverse Pinsker's inequality for general probability measures $P, Q$ as follows.
	\begin{lemma}
		If $P\ll Q$ and let $\beta_1,\beta_2\in[0,1]$ be given by
		\begin{equation}\label{reverse condition 1}
		\beta_1^{-1}:=\sup_{x\in\Omega}\frac{dP}{dQ}(x),\;\;\beta_2:=\inf_{x\in\Omega}\frac{dP}{dQ}(x).
		\end{equation}
		Then 
		\begin{equation}
		\KL(P\|Q)\le -\frac{1}{2}\left(\frac{\ln \beta_1}{1-\beta_1}-\beta_2   \right)\|P-Q\|_{TV}.
		\end{equation}
	\end{lemma}
	Sason  also discussed the discrete case in \cite{reversePin2015} that
	\begin{lemma}
		If $P, Q$ are discrete probability measures on a countably infinite set $B$, then
		\begin{equation}
		\KL(P\|Q)\le \ln\left( 1+\frac{\|P-Q\|_{TV}^2}{2 \min_{x\in B}Q(x)} \right).
		\end{equation}
	\end{lemma}

For two probability measure $P$ and $Q$ defined on a common measurable space $(\Omega, \mathcal{F})$ such that $P\ll Q$,
we say that $(P,Q)\in \mathcal{A}(\delta,m,M)$ if 
	\begin{equation}\label{reverse condition 2}
\mbox{ess}_Q\inf\frac{dP}{dQ}=m,\;\;	\mbox{ess}_Q\sup\frac{dP}{dQ}=M,\;\; \|P-Q\|_{TV}=\delta,
\end{equation}
where $\mbox{ess}_Q\inf f:=\sup\{ b\in\R| Q(\{x:f(x)<b\})=0\}, \mbox{ess}_Q\sup f:=\inf\{ b\in\R| Q(\{x:f(x)>b\})=0\}$.
 Binette(2019) \cite{Renyipinsker2019} showed a reverse Pinsker's inequality by an equality form as follows.
	\begin{lemma}
		If $\delta\ge 0, m\ge 0, M<\infty$ and $\mathcal{A}(\delta,m,M)\neq \emptyset$, then
		\begin{equation}
		\sup_{(P,Q)\in \mathcal{A}(\delta,m,M)} D_f(P\|Q)=\delta\left(\frac{f(m)}{1-m}+\frac{f(M)}{M-1}  \right),
		\end{equation}
		where $f$ is any measurable function such that $f:[0,\infty)\to (-\infty, \infty]$, $f(1)=0$ and
		\begin{equation}
		D_f(P\|Q):=\E_Q\left[ f\left( \frac{dP}{dQ}\right) \right].
		\end{equation}
	\end{lemma}
If $f(t)=t\ln t$, then $	D_f(P\|Q)=	\KL(P\|Q).$

 In addition, Berend and Kontorovich(2012) \cite{KLfromPpinsker2012} discussed the  minimum KL-divergence conditional on some distances in total variation distance sense, that is,
\begin{lemma}
	Given a probability measure $P$ on  measurable space $(\Omega, \mathcal{F})$ and define 
	\begin{equation}
	\beta:=\inf\{P(A); A\in\mathcal{F}, P(A)\ge 1/2  \}.
	\end{equation}
	Then 
	\begin{enumerate}
		\item If $\beta>1/2$, than for $ 2\beta-1>\varepsilon\to 0$,
	\begin{equation}
	\inf_{Q; \|P-Q\|_{TV}\ge \varepsilon} \KL(P\|Q)=\frac{1}{8\beta(1-\beta)}\varepsilon^2-\frac{2\beta-1}{24\beta^2(1-\beta)^2}\varepsilon^3+O(\varepsilon^4).
	\end{equation}
	
	\item If $\beta=1/2$, then for $ 1>\varepsilon\to 0$,
		\begin{equation}
	\inf_{Q; \|P-Q\|_{TV}\ge \varepsilon} \KL(P\|Q)=\varepsilon^2/2+\varepsilon^4/4+O(\varepsilon^6).
	\end{equation}
	\end{enumerate}
\end{lemma}
Thus for any $P$ with $\beta=1/2$, and any small $\varepsilon$, we can find a $Q_\varepsilon$ that
\begin{equation}
 \KL(P\|Q_\varepsilon)\le  \frac{1}{2}\|P-Q_\varepsilon\|_{TV}^2+O(\varepsilon^4).
\end{equation}

The reverse Pinsker's inequality attempts to utilize only total variation distance in order to control KL-divergence and subsequently deduce an equivalent convergence relation 
$$\KL(P_n\|Q)\to 0\Leftrightarrow \|P_n-Q\|_{TV}\to 0,\;\;\mbox{as}\;\;n\to\infty,$$
by combining it with the Pinsker's inequality. However, it is evident that the conditions of finiteness and non-zero values for the parameters \eqref{reverse condition 1} and \eqref{reverse condition 2} are strong assumptions. In this paper, we focus on random variables with densities on $\R$ and introduce the $L^2$ norm of the density to establish upper bounds on KL-divergence for a more generalized scenario. 

The paper outline and our main contributions are stated as follows. In Section \ref{section2}, we place greater emphasis on the connections between different convergence senses implied by our upper bounds  and ultimately observe that convergence in terms of KL-divergence is bounded between convergence in the $L^1$ norm  the $L^2$ norm of the density, as shown in Theorem \ref{general limit for KL from L1,2} and Theorem \ref{with 0 point}. We also establish an equivalence between the aforementioned convergence senses under certain assumptions, as Corollary \ref{general  equivalence for different convergent forms}. Particularly, we focus on a specific case that we consider the measures of difference between random variable $X$ and Gaussian random variable $G$, and deduce the above results for this case as Theorem \ref{L1,L2 bound for G} and Corollary \ref{convergence equaivalence for G}. The key Corollary \ref{KL and ||infty} are also proved in Section \ref{section2}, which is treated as the key tool to obtain a rate theorem for conditional central limit theorem. Moving on to Section \ref{section3}, our upper bounds are applied to demonstrate the   central limit theorem, with valuable assistance from the local limit theorem. Theorem \ref{KL after local limit theorem} is set for  central limit theorem with independent sum case, while after some truncation, the rate theorem for conditional central limit theorem is established as Theorem \ref{rate theorem ||infty}.

	\section{New Upper bounds for KL-divergence}\label{section2}
In this section, we present novel upper bounds for the KL-divergence $\KL(X\|Y)$ between two continuous random variables $X$ and $Y$, leveraging the $L^1$ and $L^2$ norms of their densities. We further expand these bounds when $Y$ follows a Gaussian distribution. Our bounds establish encompassing relationships between various types of convergences, including KL-divergence convergence, $L^1$ norm convergence, and $L^2$ norm convergence for continuous random sequences. In the following part, for any $1\le p\le\infty$,  any function $f$ and any measurable set $A$, we denote $L^p(\R)$ by $L^p$ and $\int_A f(x)dx$ by $\int_A f$ for simplicity.
	
		\subsection{The General Case}
	Firstly, we consider two random variables with one having non-zero density, and we have the following result.
	 
	\begin{theorem}\label{general limit for KL from L1,2}
		Given two random variables $X, Y$, with densities $p_X,p_Y$, where $0<p_Y(y)\le D<\infty,\;\;\forall y\in\R$. Then for any set $\mathcal{A}:=[-A, A]\subset \R, A>0,$ and $ s>1$,
		\begin{equation}\label{L2 as bound for KL}
		\begin{aligned}
		&\KL(X\|Y)\\
		&\le (1+|\ln D|)[\P(|Y|\ge A)+\|p_X-p_Y\|_1]\\
		&\;\;\;\;+\E_X(|\ln p_Y(X)|^s)^{1/s}[\P(|Y|\ge A)+\|p_X-p_Y\|_1]^{1-1/s}
		(1+\max_{y\in\mathcal{A}}p^{-1}_Y(y))\|p_X-p_Y\|_2^2.
		\end{aligned}
		\end{equation}
	\end{theorem}
	\pf
	Define 
	$$\mathcal{B}:=\mathcal{A}^c\cap\{x; p_X(x)\le p_Y(x)\},\;\; \mathcal{C}:=\mathcal{A}^c\cap\mathcal{B}^c.$$
	Since $t\ln t\le (t-1)+(t-1)^2, \forall t>0$, we have
	\begin{align}
	&\KL(X\|Y)=\int_\mathcal{A} p_{X}\ln\frac{p_X}{p_Y}+\int_\mathcal{B} p_X\ln\frac{p_X}{p_Y}+\int_\mathcal{C} p_X\ln\frac{p_X}{p_Y}\\
	&\le \int_{\mathcal{A}} p_Y \frac{p_X}{p_Y}\ln \frac{p_X}{p_Y}+\int_{\mathcal{C}} p_X \ln \frac{p_X}{p_Y}\\
	&\le \int_{\mathcal{A}} \left( \frac{p_X}{p_Y}-1\right)  p_Y +\int_{\mathcal{A}} \left( \frac{p_X}{p_Y}-1\right) ^2 p_Y+\int_{\mathcal{C}} p_X \ln \frac{p_X}{p_Y}\\
	&=\int_{\mathcal{A}^c} (p_Y-p_X)+\int_{\mathcal{A}}\frac{(p_X-p_Y)^2}{p_Y} +\int_{\mathcal{C}} p_X \ln \frac{p_X}{p_Y}\\
	&\le \P(|Y|\ge A)+(\max_{y\in\mathcal{A}}p^{-1}_Y(y))\|p_X-p_Y\|_2^2+\int_\mathcal{C} p_X\ln( 1+|p_X-p_Y|p_Y^{-1})\\
	&\le \P(|Y|\ge A)+(\max_{y\in\mathcal{A}}p^{-1}_Y(y))\|p_X-p_Y\|_2^2+\int_\mathcal{C} p_X|\ln p_Y|\\
	&\;\;\;\;+\int_\mathcal{C} p_X\ln( D+|p_X-p_Y|)\\
	&\le \P(|Y|\ge A)+(\max_{y\in\mathcal{A}}p^{-1}_Y(y))\|p_X-p_Y\|_2^2+\int_{x\in\mathcal{A}^c} p_X|\ln p_Y|\\
	&\;\;\;\;+\int_\mathcal{C} p_X(|\ln D|+D^{-1}|p_X-p_Y|)\\
	&\le \P(|Y|\ge A)+(\max_{y\in\mathcal{A}}p^{-1}_Y(y))\|p_X-p_Y\|_2^2+\int_{x\in\mathcal{A}^c} p_X|\ln p_Y|\\
	&\;\;\;\;+|\ln D|\|p_X-p_Y\|_1+\|p_X-p_Y\|_2^2+\int_\mathcal{C} p_Y(|\ln D|+D^{-1}|p_X-p_Y|)\\
	&\le (1+|\ln D|)\P(|Y|\ge A)+(1+\max_{y\in\mathcal{A}}p^{-1}_Y(y))\|p_X-p_Y\|_2^2+\int_{x\in\mathcal{A}^c} p_X|\ln p_Y|\\
	&\;\;\;\;+|\ln D|\|p_X-p_Y\|_1+\int_\mathcal{C} D(D^{-1}|p_X-p_Y|)\\
	&\le (1+|\ln D|)\P(|Y|\ge A)+(1+\max_{y\in\mathcal{A}}p^{-1}_Y(y))\|p_X-p_Y\|_2^2\\
	&\;\;\;\;+\int_{x\in\mathcal{A}^c} p_X|\ln p_Y|+(|\ln D|+1)\|p_X-p_Y\|_1.
	\end{align}
	The  inequality now comes from 
	\begin{equation}
	\begin{aligned}
	\int_{x\in\mathcal{A}^c} p_X|\ln p_Y|&\le \E_X(|\ln p_Y(X)|^s)^{1/s}\P^{1-1/s}(|X|\ge A)\\
	&\le \E_X(|\ln p_Y(X)|^s)^{1/s}[\P(|Y|\ge A)+\|p_X-p_Y\|_1]^{1-1/s}.
	\end{aligned}
	\end{equation}
	 
	\e

	Directly, we deduce that $\|\cdot\|_2$ can controls KL-divergence, as the following corollary.
	\begin{corollary}\label{c-2.1}
		Given a random variable $ Y$ with density $p_Y$, where $0<p_Y(y)\le D<\infty,\;\;\forall y\in\R$, and a random variable sequence $\{X_n, n\ge 1\}$ with densities $\{p_n, n\ge 1\}$. If $\|p_n-p_Y\|_2\to0$ and there exists $s>1$ such that $T_n:=|\ln p_Y(X_n)|$ has  uniformly bounded $s$-th moment, then $\KL(X_n\|Y)\to 0.$
	\end{corollary}
	\pf
	We use Theorem \ref{general limit for KL from L1,2} for all $X_n$. 
	Since $\|p_n-p_Y\|_2\to0$, we can always choose $\mathcal{A}_n=[-A_n, A_n]$ s.t. 
	$$\max_{y\in[-A_n,A_n]}p_Y^{-1}(y) \le \min\{M, \|p_n-p_Y\|^{-1}_2 \}$$
	with $\lim_{n\to\infty}A_n=\infty,$ where
	$$M:=\max_{y\in[-1,1]}p_Y^{-1}(y).$$
	Then with above $\mathcal{A}_n$, the assumptions imply that as $n\to\infty$,
	\begin{align*}
	a_n:=&(1+e^{-1}+|\ln D|)\P(|Y|\ge A_n)+(1+\max_{y\in\mathcal{A}_n}p^{-1}_Y(y))\|p_{n}-p_Y\|_2^2\\
	&+\E_{X_n}(|\ln p_Y(X_n)|^s)^{1/s}\P^{1-1/s}(|Y|\ge A_n)
	\end{align*}
tends to 0.  Further, we have $\|p_n-p_Y\|_2\to 0 \Rightarrow \|p_n-p_Y\|_1\to 0$,  which is proved in \ref{A1}, and thus the proof completes.
	\e
	
	\begin{remark}
The assumption ``$\E T^s_n<C, \forall n"$ is not hard to be satisfied, for example, it holds when $X_n$s are Gaussian random variables with bounded variances.
	\end{remark}
	
%
%
	
	Since $\|\cdot\|_2$ control KL-divergence, and KL-divergence control $\|\cdot\|_1$, once the convergence of  $\|\cdot\|_1$ can deduce the convergence of $\|\cdot\|_2$, the three convergent become equivalent. Specifically, we deduce the following equivalence between different convergent forms.

	\begin{corollary}\label{general  equivalence for different convergent forms}
		Given a random variable $ Y$ with density $p_Y$, where $0<p_Y(y)\le D<\infty,\;\;\forall y\in\R$, and a random variable sequence $\{X_n, n\ge 1\}$ with densities $\{p_n, n\ge 1\}$. If $\|p_n-p_Y\|_2\to0$ and there exists $s>1$ such that $T_n:=|\ln p_Y(X_n)|$ has  uniformly bounded $s$-th moment, and one of the following assumptions holds:
		\begin{enumerate}
				\item $\exists 2<q\le \infty, \sup_n\|p_n-p_Y\|_q<\infty.$
			
			\item $\limsup_n \|p_n-\|_2\le \|p_Y\|_2<\infty$.
		
		\end{enumerate}
		Then we obtain the equivalence of the following statements:
		\begin{enumerate}
			\item $p_n\xrightarrow{L^1(\R)} p_Y,\;\; n\to\infty.$
			
			\item $p_n\xrightarrow{L^2(\R)} p_Y,\;\; n\to\infty.$
			
			\item $\KL(X_n\|Y)\to 0,\;\; n\to\infty.$
		\end{enumerate}
		
	\end{corollary}
		\pf Due to Pinsker's inequality and Corollary \ref{c-2.1}, we only need to prove ``$p_n\xrightarrow{L^1(\R)} p_Y, n\to\infty"$ $\Rightarrow$ ``$p_n\xrightarrow{L^2(\R)} p_Y, n\to\infty"$ under the assumptions. 
			\begin{enumerate}
			\item  If $\exists 2<q\le \infty, \sup_n\|p_n-p_Y\|_q<\infty,$ then we can find $\alpha\in(0,1]$ \textit{s.t.} $\alpha+\frac{1-\alpha}{q}=\frac{1}{2}$. Then by H$\ddot{\mbox{o}}$lder's inequality, 
			\begin{equation}
			\|p_n-p_Y\|_2\le \|p_n-p_Y\|_1^\alpha\|p_n-p_Y\|_q^{1-\alpha}\to 0,\;\;n\to\infty.
			\end{equation}

		\item  If $\limsup_n \|p_n\|_2\le \|p_Y\|_2<\infty$, 	we know that $(p_n)$ is bounded in $L^2$. According to \cite[Exercise 4.16.2]{functionanalysis} (which has been proved in \ref{A2}), we have $p_n\rightharpoonup p_Y$ weakly $\sigma(L^2,L^2)$. By \cite[Proposition 3.5]{functionanalysis}, we deduce $\|p_n\|_2\to \|p_Y\|_2, n\to\infty$.   By \cite[Exercise 4.19.1]{functionanalysis} (which has been proved in \ref{A2}), we finish the proof.
		\end{enumerate}

	\e

	Further more,	if $Y$ has some zero density points on $\R$ but $X\ll Y$, we have a similar result as Theorem \ref{general limit for KL from L1,2}.
		\begin{theorem}\label{with 0 point}
			Given $X, Y$, with density $p_X,p_Y$, where $p_Y(y)\le C<\infty,\;\;\forall y\in\R$ and $X\ll Y$. Then for any set $\hat{\mathcal{A}}:=\{y;p_Y(y)\ge \frac{1}{A}\}\subset \R, A>0,$ and $ s>1$,
			\begin{equation}
			\begin{aligned}
			&\KL(X\|Y)\\
			&\le (1+|\ln C|)[\P(Y\in \hat{\mathcal{A}}^c)+\|p_X-p_Y\|_1]\\
			&\;\;\;\;+\E_X(|\ln p_Y(X)|^s)^{1/s}[\P(Y\in\hat{\mathcal{A}}^c)+\|p_X-p_Y\|_1]^{1-1/s}(1+A)\|p_X-p_Y\|_2^2.
			\end{aligned}
			\end{equation}
		\end{theorem}
		Apparently, the above Corollary \ref{c-2.1} and Corollary \ref{general  equivalence for different convergent forms} can also be obtained under these conditions.

	\subsection{A Special Case}
	In this part, we consider a special case that $Y$ in Theorem \ref{general limit for KL from L1,2} is Gaussian with the same mean and variance of $X$. We denote a  Gaussian random variable with the same mean and variance of $X$ by $G_X$. 
	Under this special case, we write $\KL(X\|G_X)$ as $D(X)$ for simplify, similarly with \cite{Esseenbound}. Since for any $a\in\R, b>0$, and random variable $X$,
	\begin{equation}\label{D(ax)}
	D(X)=D(b(X-a)),
	\end{equation}
	we first focus on  $X$ be a random variable with density $p$ and $\E X=0, \V(X)=1$. Denote standard Gaussian random variable by $G$ and its density by $\phi(x)=(\sqrt{2\pi})^{-1}e^{-x^2/2}.$

	We   define three sets for a given $A>0$:
	$$\mathcal{A}:=[-A,A],\;\; \mathcal{B}:=\mathcal{A}^c\cap\{x; p(x)\le \phi(x)\},\;\; \mathcal{C}:=\mathcal{A}^c\cap\mathcal{B}^c.$$
	On $\mathcal{A}=[-A,A]$, we use the following lemma.
	\begin{lemma}[Lemma 2.2 in \cite{Esseenbound}]
		 $$\int_\mathcal{A} p\ln\frac{p}{\phi}\le \P(|G|\ge A)+\sqrt{2\pi}\int_{\mathcal{A}} (p-\phi)^2e^{x^2/2}-\int_{\mathcal{A}^c} p.$$
	\end{lemma}
	Then we deduce that
	\begin{equation}\label{21}
	\int_\mathcal{A} p\ln\frac{p}{\phi}\le \P(|G|\ge A)+\sqrt{2\pi}e^{A^2/2}\|p-\phi\|^2_2.
	\end{equation}
	 On $\mathcal{A}^c$, when $p/\phi\le 1$(\textit{i.e.} on $\mathcal{B}$), we have $|\ln(p/\phi)|\le e^{-1}(\phi/p).$ Finally on $\mathcal{C}$, since
	\begin{equation}\label{22}
	\frac{p(x)}{\phi(x)}\le 1+\|p-\phi\|_\infty \phi(x)^{-1}= 1+\sqrt{2\pi}\|p-\phi\|_\infty e^{x^2/2}\le \sqrt{2\pi}(\|p-\phi\|_\infty+1) e^{x^2/2},
	\end{equation}
	we can show that
	\begin{equation}\label{use infty norm}
	\int_{\mathcal{C}}p\ln(p/\phi)\le \ln(\sqrt{2\pi}(\|p-\phi\|_\infty+1))\P(X\in\mathcal{A}^c)+\frac{1}{2}\E(X^21_{X\in\mathcal{A}^c}).
	\end{equation}
	Combining  \eqref{21} \eqref{22}, similar to the proof of Theorem \ref{general limit for KL from L1,2},  we obtain the following  theorem.
		\begin{theorem}\label{L1,L2 bound for G}
		Given a random variable $X$ with density $p$, $\E X=0, \V(X)=1$ and $\E|X|^s<\infty$ for some $s>2$, then $\exists C_1,C_2,C_3>0$, such that for any $A>0$,
		\begin{equation}
		\begin{aligned}
		D(X)\le &C_1e^{A^2/2}\|p-\phi\|^2_2+(C_2+\|p-\phi\|_\infty)\left( \|p-\phi\|_1+e^{-A^2/2} \right)\\
		&+C_3(\E|X|^s)^{2/s}\left( \|p-\phi\|_1+e^{-A^2/2} \right)^{1-2/s}.
		\end{aligned}
		\end{equation}
	\end{theorem}
\pf
	\begin{align}
	&D(X)\\
	&=\int_\mathcal{A} p\ln\frac{p}{\phi}+\int_\mathcal{B} p\ln\frac{p}{\phi}+\int_\mathcal{C} p\ln\frac{p}{\phi}\\
	&\le \P(|G|\ge A)+\sqrt{2\pi}e^{A^2/2}\|p-\phi\|^2_2\\
	&\;\;\;\;+\ln(\sqrt{2\pi}(\|p-\phi\|_\infty+1))\P(X\in\mathcal{A}^c)+\frac{1}{2}\E(X^21_{X\in\mathcal{A}^c})\\
	&\le \P(|G|\ge A)+\sqrt{2\pi}e^{A^2/2}\|p-\phi\|^2_2\\
	&\;\;\;\;+\ln(\sqrt{2\pi}(\|p-\phi\|_\infty+1))\P(X\in\mathcal{A}^c)+\frac{1}{2}\E(X^21_{X\in\mathcal{A}^c})\\
	&\le e^{-A^2/2}+C_1e^{A^2/2}\|p-\phi\|^2_2+(C_2+\ln(\|p-\phi\|_\infty+1))\P(X\in\mathcal{A}^c)\\
	&\;\;\;\;+C_3(\E|X|^s)^{2/s}\P^{1-2/s}(X\in\mathcal{A}^c)\\
	&\le e^{-A^2/2}+C_1e^{A^2/2}\|p-\phi\|^2_2+(C_2+\|p-\phi\|_\infty)\P(X\in\mathcal{A}^c)\\
	&\;\;\;\;+C_3(\E|X|^s)^{2/s}\P^{1-2/s}(X\in\mathcal{A}^c).
	\end{align}
	Since
	\begin{equation}\label{prob for |x|>A, weak}
	\P(X\in\mathcal{A}^c)=\int_{|x|\ge A} p=\int_{|x|\ge A} p-\phi+\int_{|x|\ge A} \phi\le \|p-\phi\|_1+e^{-A^2/2},
	\end{equation}
	the proof completes.
	\e

	By taking $e^{A^2/2}=\|p-\phi\|^{-1}_2 $ when $\|p-\phi\|_2\le 1$,  the following corollary holds.
	\begin{corollary}
		Given a random variable $X$ with density $p$, $\E X=0, \V(X)=1$ and $\E|X|^s<\infty$ for some $s>2$, and 
		\begin{equation}\label{<1}
		\|p-\phi\|_2\le 1.
		\end{equation}
		Then $\exists C_1, C_2>0$ such that
		\begin{equation}\label{23}
		\begin{aligned}
		&D(X)\\
		&\le (C_1+\|p-\phi\|_\infty)\left( \|p-\phi\|_1+\|p-\phi\|_2 \right)+C_2(\E|X|^s)^{2/s}\left( \|p-\phi\|_1+\|p-\phi\|_2 \right)^{1-2/s}.
		\end{aligned}
		\end{equation}
	\end{corollary}
\begin{remark}
	The constant 1 in \eqref{<1} is only for simplicity, in fact for all the $X$ satisfying $\|p-\phi\|_2\le C_0$, by taking $e^{\frac{A^2}{2}}=\left(\frac{\|p-\phi\|_2}{2C_0}\right)^{-1}$, \eqref{23} still holds.
\end{remark}
	Moreover, if we want to throw away $\|p-\phi\|_\infty$, we can replace \eqref{use infty norm} by
	\begin{equation}
	\begin{aligned}
	\int_{\mathcal{C}}p\ln(p/\phi)&\le \ln(\sqrt{2\pi})\P(X\in\mathcal{A}^c)+\int_\mathcal{C}\ln(p-\phi+1)p+\frac{1}{2}\E(X^21_{X\in\mathcal{A}^c})\\
	&\le \ln(\sqrt{2\pi})\P(X\in\mathcal{A}^c)+\int_\mathcal{C}(p-\phi)p+\frac{1}{2}\E(X^21_{X\in\mathcal{A}^c})\\
	&\le \ln(\sqrt{2\pi})\P(X\in\mathcal{A}^c)+\|p-\phi\|_2^2+\int_{\mathcal{C}}(p-\phi)\phi+\frac{1}{2}\E(X^21_{X\in\mathcal{A}^c})\\
	&\le \ln(\sqrt{2\pi})\P(X\in\mathcal{A}^c)+\|p-\phi\|_2^2+\|p-\phi\|_1e^{-A^2/2}+\frac{1}{2}\E(X^21_{X\in\mathcal{A}^c}).
	\end{aligned}
	\end{equation}
	We thus obtain this result.
	\begin{corollary}
		Given a random variable $X$ with density $p$, $\E X=0, \V(X)=1$ and $\E|X|^s<\infty$ for some $s>2$, and 
		$$\|p-\phi\|_2\le 1.$$
		Then $\exists C_1, C_2>0$ such that
		\begin{equation}
		D(X)\le C_1\left( \|p-\phi\|_1+\|p-\phi\|_2 \right)+C_2(\E|X|^s)^{2/s}\left( \|p-\phi\|_1+\|p-\phi\|_2 \right)^{1-2/s}.
		\end{equation}
	\end{corollary}
\begin{remark}
	The constant 1   is also only for simplicity.
\end{remark}
	
Similar to Corollary \ref{general  equivalence for different convergent forms}, we have the following result between different convergences.

	\begin{corollary}\label{convergence equaivalence for G}
		Given a sequence of  random variables $\{X_n, n\ge 1\}$ with densities $\{p_n, n\ge 1\}$, $\E X_n=0, \V(X_n)=1, \forall n\ge 1$ and $\exists s>2, \sup_n\E|X_n|^s<\infty $, then if one of the following assumptions holds:
		\begin{enumerate}
	\item $\exists 2<q\le \infty, \sup_n\|p_n-p_Y\|_q<\infty.$
	
	\item $\limsup_n \|p_n\|_2\le \|p_Y\|_2<\infty$.
	
\end{enumerate}
	We can obtain the equivalence of the following statements:
		\begin{enumerate}
			\item $p_n\xrightarrow{L^1(\R)} \phi,\;\; n\to\infty.$
			
			\item $p_n\xrightarrow{L^2(\R)} \phi,\;\; n\to\infty.$
			
			\item $\KL(X_n\|G)\to 0,\;\; n\to\infty.$
		\end{enumerate}
	\end{corollary}

	When $\|p-\phi\|_\infty\le 1/2$, we have
	\begin{equation}
	\|p-\phi\|_1=2\int (\phi-p)^+ \le 4A\|p-\phi\|_\infty+2\frac{1}{A}e^{-A^2/2}\le C\sqrt{|\ln(\|p-\phi\|_\infty)|}\|p-\phi\|_\infty,
	\end{equation}
	where we take $A=\sqrt{2|\ln(\|p-\phi\|_\infty)|}$, then
	\begin{equation}
	\|p-\phi\|_2\le \|p-\phi\|_\infty\|p-\phi\|_1\le C\sqrt{|\ln(\|p-\phi\|_\infty)|}\|p-\phi\|_\infty^2.
	\end{equation}
	Thus we conclude the following corollary which is the key to the next section.
	\begin{corollary}\label{KL and ||infty}
		Given a random variable $X$ with density $p$ satisfying $\E X=0, \V(X)=1$ and $\E|X|^s<\infty$ for some $s>2$, and  $\|p-\phi\|_\infty\le 1/2$, then $\exists C_1>0, \forall \varepsilon>0$,
		\begin{equation}
		\begin{aligned}
				D(X)&\le C_1(\E|X|^s)^{2/s}\left(  \sqrt{|\ln(\|p-\phi\|_\infty)|}\|p-\phi\|_\infty\right)^{1-2/s}\\
				&\le C_1M(\varepsilon)(\E|X|^s)^{2/s}\left(  \|p-\phi\|_\infty\right)^{(1-2/s)(1-\varepsilon)},
		\end{aligned}
		\end{equation}
		where
		$$M(\varepsilon):=\max_{x\in[0,1]} x^\varepsilon \sqrt{|\ln x|}+1<\infty.$$
	\end{corollary}
	Thanks to Corollary \ref{KL and ||infty}, we can prove the rate theorem for conditional central limit theorem, see Theorem \ref{rate theorem ||infty} in Section \ref{section3}.

	\section{Conditional Central Limit Theorem}\label{section3}

As emphasized in \cite{yuan2014}, the random nature of many problems arising in the applied sciences leads to mathematical models where conditioning is
present, thus it is necessary to study the conditional distribution.

Holst (1979) proved a conditional central limit theorem (CCLT) for a special conditional distribution in \cite{Holst1979TwoCL}, which has a simple form that under some regular assumptions for $(\xi_i, \eta_i)$,
$$\P\left( \sum_{i=1}^n \frac{\xi_i}{\sqrt{\sum_{i=1}^n\V(\xi_i)}}\le t{\bigg|}\sum_{i=1}^n \eta_i=n\E \eta_1+o(1)\right) \to \Phi(t), \;\; n\to\infty,$$
where $\eta_i$'s are \textit{i.i.d.} discrete random variables, and $\xi_i$'s are also \textit{i.i.d.} with $E \xi_i=0,\forall i$, $\Phi$ is standard Gaussian distribution function. The result was strengthened in \cite{Janson2001MomentCI} to the case  $n\E  \eta_1+O(\sqrt{n}).$

In \cite{Rubshtein1996} Rubshtein proved a CCLT for ergodic stationary process. Using martingale method, Rubshtein showed that for a two-dimensional ergodic stationary process $(\bm{\xi},\bm{\eta}):=\{(\xi_n,\eta_n),n\geq  1\}$ satisfying some regular conditions, the conditional distribution function
\begin{equation}
F_{n,\bm{z}}(t) = \P\left(\frac{S_n - \E (S_n|\bm{\eta}=\bm{z})}{\sqrt{\V (S_n|\bm{\eta}=\bm{z})}}\leq t| \bm{\eta}=\bm{z}\right), t\in\R
\end{equation}
converges weakly to the distribution function of standard Gaussian for realizations $\bm{z} =\bm{\eta}(\omega)= \{\eta_n(\omega)\}$ almost surely, where $S_n = \xi_1+\xi_2+\cdots+\xi_n$.

More recently, \cite{yuan2014} showed a limit theorem for the conditionally independent case. \cite{2021Stein} used Stein method to deduce an upper bound for Wasserstein distance between conditional distribution and Gaussian distribution. 

It is worthy to note that Yao \textit{et al.} (2024) use entropy method to prove a conditional central limit theorem in \cite{ECCLT} which covers the independent sum case. Yao \textit{et al.} also deduce a rate theorem for conditional central limit theorem in \cite{ECCLT} according to the Berry-Esseen bound in \cite{Esseenbound}. However, a finite 64-th moment condition is needed (see details in Remark \ref{64} and \cite{ECCLT}), which is necessary to be further improved.

	In this section, we use our upper bounds in Section \ref{section2} to deduce a conditional central limit theorem in KL-divergence sense (see Theorem \ref{KL after local limit theorem}), based on local limit theorem stated in  the next part, and eventually prove a new rate theorem for conditional central limit theorem (see Theorem \ref{rate theorem ||infty}).

	\subsection{Local Limit Theorem}
  We ask the help for traditional local limit theorem, specifically, 	we combine the classical results in \cite{Sumofindependent}: Theorem 10 in Chapter VI, Lemma 1 in Chapter V and Lemma 10 in Chapter VI, and then conclude the following lemma.
	\begin{lemma}\label{||infty}
		Consider an independent sequence  of random variables $\{X_n, n\ge 1\}$ with densities $\{p_{X_n}, n\ge 1\}$, and $\E X_n=0, \forall n\ge 1$. Denote
		\begin{equation}
		B_n:=\sum_{i=1}^n \V(X_i).
		\end{equation}
		If the following conditions hold:
		\begin{enumerate}
			\item $\lim_{n\to\infty}B_n=\infty$. 
			
			\item $\sum_{i=1}^n \E |X_i|^3=O(B_n).$
			
			\item There exists a subsequence $X_{n_1}, X_{n_2}\cdots, X_{n_m},\cdots$, \textit{s.t.} for some $\lambda>0$,
			\begin{equation}
			\mu:=\liminf_{n\to\infty} \frac{\#\{m:n_m\le n  \}}{n^\lambda}>0,
			\end{equation}
			and $\sup_{m} \|p_{X_{n_m}}\|_{TV}<\infty$, where $\#$ denotes the counting measure and $\|\cdot\|_{TV}$ denotes the total variation of a function, \textit{i.e.}
			$$\|f\|_{TV}:=\sup_{n, a_1<a_2<\cdots<a_n}\sum_{i=1}^n|f(a_i)-f(a_{i+1})|.$$
		\end{enumerate}
		
		Then if we denote the density function of $S_n:=\sum_{i=1}^n X_i$ by $p_n$, we obtain
		\begin{equation}
		\|p_{n}-\phi\|_\infty \le  C_0B_n^{-3/2}\sum_{i=1}^n \E|X_i|^3+C_1 \sqrt{B_n} e^{-C_2n^\lambda},
		\end{equation}
		where $\phi=(\sqrt{2\pi})^{-1}e^{-x^2/2}$ is the density of standard normal distribution, $C_0, C_1, C_2$ are positive constants.
	\end{lemma}
\pf
Theorem 10 in Chapter VI in \cite{Sumofindependent} says that
\begin{equation}
\|p_n-\phi\|_\infty \le C_0B_n^{-3/2}\sum_{i=1}^n \E|X_i|^3+ \sqrt{B_n}\int_{|t|>\frac{1}{2}K_1} \prod_{j=1}^n |f_n(t)|dt,
\end{equation}
where $C_0$ is an absolutely constant, $K_1:=\sup_n B_n^{-1}\sum_{i=1}^n \E|X_i|^3$ and $f_n$ is the characteristic function of $X_n$. Lemma 10 in Chapter VI of \cite{Sumofindependent} says that if condition 3 holds, there exists a constant $C_B$ such that 
\begin{equation}\label{C_B}
\int_{|t|>\frac{1}{2}K_1} \prod_{j=1}^n |f_n(t)|dt\le K_2\frac{1}{\sqrt{\gamma \mu n^\lambda}}e^{-\gamma \mu n^\lambda \frac{1}{4}K_1}+K_3\left( \frac{1}{2}\right) ^{\mu n^{\lambda}}4C_B^2,
\end{equation}
where $K_2, K_3$ are absolutely constants, and $\gamma=\frac{3}{128C_B^2}$  according to Theorem 1 in Chapter I of \cite{Sumofindependent}. Clearly, proof completes immediately.
\e

For a random variable $X$  absolutely continuous density function $p$, the Fisher information defined as follows is able to control its total variance.
\begin{definition}[Fisher information]
	For a random variable $X$  absolutely continuous density function $p$, its Fisher information is defined as 
	\begin{equation}
	J(X)=J(p):=\int_\R \frac{|p'|^2}{p}.
	\end{equation}
\end{definition}
In fact, for any absolutely continuous density function $p$, we always have
		\begin{equation}\label{TV and c.h.f}
		\|p\|_{TV}\le \int_\R |p'|=\int_\R \sqrt{p}\frac{|p'|}{\sqrt{p}}  \le 1\cdot \sqrt{\int_\R \frac{|p'|^2}{p}}=J(p),
		\end{equation}
		\textit{i.e.} finite Fisher information can deduce the bounded variation condition in Lemma \ref{||infty}.

	Next we discuss   the case that $\sup_n \|p_{X_n}\|_{TV}=\infty$, \textit{i.e.} $C_B=C_B(n)$ is not a constant. Indeed, given a random variable $X$ with absolutely continuous density $p$ and  characteristic function $f$,   we use integration by parts to obtain
	\begin{align}
	|f(t)|=|\int_\R p(x)e^{itx}dx|=|\frac{1}{it}\int_\R p(x)de^{itx}|=|\frac{1}{it}\int_\R e^{itx}p'(x)dx|\le \frac{1}{t}\int_\R |p'|.
	\end{align}
	Thus $C_B(n)$ can be chosen as $\sup_n \int_\R |p_{X_n}'|$ or $\sup_n \sqrt{J(X_n)}$.

%

	Recall \eqref{C_B},  we obtain the following proposition for the case $C_B(n)$ is not a constant.
	\begin{theorem}\label{constant estimation for C(J)}
	Consider an independent sequence  of random variables $\{X_n, n\ge 1\}$ with absolutely continuous densities $\{p_{X_n}, n\ge 1\}$, and $\E X_n=0, \forall n\ge 1$. Denote
\begin{equation}
B_n:=\sum_{i=1}^n \V(X_i).
\end{equation}
If the following conditions hold:
\begin{enumerate}
	\item $\lim_{n\to\infty}B_n=\infty$. 
	
	\item $\sum_{i=1}^n \E |X_i|^3=O(B_n).$
	
	\item There exists a subsequence $X_{n_1}, X_{n_2}\cdots, X_{n_m},\cdots$, \textit{s.t.} for some $\lambda>0$,
	\begin{equation}
	\mu:=\liminf_{n\to\infty} \frac{\#\{m:n_m\le n  \}}{n^\lambda}>0,
	\end{equation}
	and 
	\begin{equation}
	\inf_m \frac{n^s}{\int_\R |p'_{X_{n_m}}|} >0
	\end{equation}
	for some $s\in(0,\lambda/2)$.
\end{enumerate}
Then if we denote the density function of $S_n:=\sum_{i=1}^n X_i$ by $p_n$, we obtain
\begin{equation}
\|p_{n}-\phi\|_\infty \le  C_0B_n^{-3/2}\sum_{i=1}^n \E|X_i|^3+C_1 \sqrt{B_n} e^{-C_2n^{\lambda-2s}},
\end{equation}
where $\phi=(\sqrt{2\pi})^{-1}e^{-x^2/2}$ is the density of standard normal distribution, $C_0, C_1, C_2$ are positive constants.
Further more, if $\lim_{n\to\infty} \frac{B_n}{n^\alpha}=0$ for some $\alpha>0$, then
\begin{equation}
\|p_{n}-\phi\|_\infty \le  C B_n^{-3/2}\sum_{i=1}^n \E|X_i|^3,
\end{equation}
where $C$ is a constant.
	\end{theorem}
\begin{remark}
	Note that the condition $\lim_{n\to\infty} \frac{B_n}{n^\alpha}=0$ holds for $\alpha=2$ when $\sup_n \E |X_n|^k<\infty$ for some $k\ge 2$.
\end{remark}
In the remain sections we only consider the case $\lambda=\mu=1$ for simplicity, and our results also hold once condition 3 in Theorem \ref{constant estimation for C(J)} is satisfied.

\subsection{Convergence rate for $\E D(S_n|y)$}
We now show the CCLT for KL-divergence sense, as an application of our bounds in Section \ref{section2}. However, before showing our CCLT, we will present a non-conditional version as a prelude firstly. Specifically, 
we use Corollary \ref{KL and ||infty} and  Theorem \ref{constant estimation for C(J)} to deduce the following theorem immediately.

\begin{theorem}\label{KL after local limit theorem}
	Consider an independent sequence  of random variables $\{X_n, n\ge 1\}$ with densities $\{p_{X_n}, n\ge 1\}$, $\E X_i=0$, $\E|X_i|^k\le M, k\ge3$ and $(\int_\R |p_{X_i}'|)^2<J,\forall i$. Denote
	$$S_n=\frac{X_1+X_2+\cdots+X_n}{\sqrt{B_n}}.,$$
	where $B_n=\sum_{i=1}^n \V(X_i)$, then for any $\varepsilon\in(0,1/2)$,
	\begin{equation}\label{upper bound for D under kth momoent}
	D(S_n)\le C_1(\varepsilon)\left(\frac{\E|\sum_{i=1}^n X_i|^k}{(B_n)^{\frac{k}{2}}}\right)^{\frac{2}{k}}\left(C\frac{\sum_{i=1}^n \E|X_i|^3}{(B_n)^{\frac{3}{2}}}\right)^{1-\frac{2}{k}-\varepsilon},
	\end{equation}
	where $C_1(\varepsilon)$ only depends on $\varepsilon$, and $C=C(J,M)$ only depends on $J, M$.
	
	In particular, if $\E|X_i|^4\le M$, then for any $\varepsilon\in(0,1/2)$,
	\begin{equation}\label{upper bound for D under 4th momoent}
	D(S_n)\le C_1(\varepsilon)\left(C\frac{\sum_{i=1}^n \E|X_i|^3}{(B_n)^{\frac{3}{2}}}\right)^{\frac{1}{2}-\varepsilon}\sqrt{\frac{\sum_{i=1}^n \E|X_i|^4}{(B_n)^{2}}+3}.
	\end{equation}
\end{theorem}


According to \eqref{TV and c.h.f} and Theorem \ref{constant estimation for C(J)},
if for any $n$, $p_{X_n}$ is absolutely continuous and $ J(p_{X_n})< n^s$ for some $s\in(0,1)$,
then the constant $C$ in \eqref{upper bound for D under kth momoent}  in Theorem \ref{KL after local limit theorem} only depends on $s, M$. Specifically, we can find $r$ such that $\sup_n M^{2/k}n^{-r}\le1$ and take
$C=C(s,r)$ depending only on $s$ and $r$.
Further more, using the CR bound
$$\V(X)\ge J(X)^{-1}$$
and Jensen's inequality, we can obtain the following corollary.
\begin{corollary}\label{certain estimation for D}
	Let $X_1,X_2,\cdots,X_n$ be an independent sequence  with absolutely continuous densities, $\E|X_i|^k\le M, k\ge 4$ and $J(X_i)\le\kappa n^s, \kappa>0, s\in(0,1)$. Then for any $\varepsilon\in(0,1/2)$, and $r>0$ 
	\begin{equation}
	D(S_n)\le C(\varepsilon,s,r,\kappa)\left( \frac{M^{4/k}}{n^{1-2s}}+3 \right)^{1/2}\left( \frac{M^{3/k}}{n^{1/2-3s/2}}  \right)^{1/2-\varepsilon},
	\end{equation}
  when $n>M^{\frac{2}{kr}}$.
\end{corollary}
To prove the rate of convergence for conditional central limit theorem, we still need the following lemma to control $D$ when $J(X)$ is very large.
\begin{lemma}
	Let $X$ be a random variable with zero mean, finite Fisher information and finite second moment, then
	\begin{equation}\label{bounds for h}
	\frac{1}{2}\ln(2\pi e\frac{1}{J(X)})\le h(X)\le \frac{1}{2}\ln(2\pi e\V(X)),
	\end{equation}
	where $h(X):=-\int_\R p_X(x)\ln p_X(x)dx$, $p_X$ is the density of $X$.
	Further more, $\forall t>0$,
	\begin{equation}\label{bound for D when J large}
	D(X)\le C(t)(\V^t(X)+J^t(X)),
	\end{equation}
	where $C(t)$ only depends on $t$.
\end{lemma}
\pf
The LHS of \eqref{bounds for h} comes from the inequality
$$\frac{1}{2\pi e}e^{2h(X)}J(X)\ge 1$$
in \cite{Thomas2018}, while the RHS of \eqref{bounds for h} is true since Gaussian random variable has the maximum differential entropy when variance is fixed.

For \eqref{bound for D when J large}, we use \eqref{bounds for h} to deduce
\begin{align}
D(X)&=\frac{1}{2}\ln(2\pi e\V(X))-h(X)\\
&\le \frac{1}{2}\ln(2\pi e\V(X))-\frac{1}{2}\ln(2\pi e\frac{1}{J(X)})\\
&=\frac{1}{2}\ln(\V(X)J(X))\\
&\overset{a}{\le} \begin{cases}
C(t)J^t(X),  &\V(X)\le 1;\\
C(t)\V^t(X)+C_1(t)J^t(X), &\V(X)>1,
\end{cases}\\
&\le C(t)\V^t(X)+C(t)J^t(X),
\end{align}
where 
\begin{equation}\label{C(t)}
C(t):=\sup_{x\ge 1} \frac{1/2\ln x}{x^t},
\end{equation}
and '$a$' holds since $\V(X)J(X)\ge 1$.
\e

For conditional distribution, we use notation $X|Y=y$ to represent the distribution of $X$ condition on the realization $Y=y$. We also use $J(X|Y=y), \V(X|Y=y), D(X|Y=y)$ to denote the related functions for conditional random variable.
Based on the above preparations, our CCLT is stated as follows.
\begin{theorem}\label{rate theorem ||infty}
	Let $(X,Y)$ be a two-dimensional random variable. Suppose $(X_1,Y_1),(X_2,Y_2),\cdots,(X_n,Y_n)$ are independent and identically distributed as $(X,Y)$. Let $S_n = (X_1 + X_2 +\cdots+X_n) / \sqrt{n}$, and $\bm{Y}_n = (Y_1,Y_2,\cdots,Y_n)$. Assume the following conditions hold.
	\begin{enumerate}
		\item $(X|Y) \prec \mathcal{L}, \ \P_Y-a.s.$, where $\mathcal{L}$ denotes Lebesgue measure and $\mu \prec \lambda$ means $\mu$ is absolutely continuous with respect to $\lambda$. Further, for $\P_Y-a.s.$, the conditional random variable $(X|Y)$ has absolutely continuous density.
		\item  (Finite conditional Fisher information and moment) $\exists u>3, k>6, \E|X|^k<\infty, \E J^u(X|Y)<\infty$, and
		\begin{equation}
		\alpha:=\frac{6}{k}+\frac{3}{u}<1,
		\end{equation}
	where
	\begin{equation}
	\E J^u(X|Y):=\int J^u(X|Y=y)\P_Y(y).
	\end{equation}
		
	\end{enumerate}
	Then for any $\varepsilon\in (0,\frac{1-\alpha}{4+\alpha})$,
	\begin{equation}\label{rate result}
	\E D(S_n|\bm{Y}_n):=\int D(S_n|\bm{Y}_n=\bm{y}_n)\P_{\bm{Y}_n}(\bm{y}_n) =O\left(\left( \frac{1}{n}\right)^{\frac{1-\alpha}{4+\alpha}-\varepsilon}\right).
	\end{equation}
\end{theorem}
\pf  Fix $\varepsilon\in(0,1/2)$. For any realization $\bm{Y}_n=\bm{y}_n=(y_1,y_2,\cdots,y_n)$, we can write the conditional distribution as
\begin{equation}
(S_n|\bm{Y}_n=\bm{y}_n)\overset{d}{=}\eta_1+\eta_2+\cdots+\eta_{n},
\end{equation}
where $\eta_i$s are independent and 
$$\eta_i\overset{d}{=}\frac{(X|Y=y_i)}{\sqrt{n}}.$$

Fix $s\in(1/u, 1)$ and define
$$A_1(n):=\{ y:   J(X|Y=y)\le n^s \},$$
$$A_2(M):=\{ y: \E(|X|^k|y)\le M\},$$
$$A(M):=(A_1(n)\bigcap A_2(M))^n=\{\bm{y_n}:y_i\in A_1(n)\cap A_2(M),\forall i \},$$
and 
then by Corollary \ref{certain estimation for D} and \eqref{D(ax)},
for any $\bm{y}_n=(y_1,\cdots,y_n)\in A(M), r\ge0$,
\begin{equation}\label{bound in AM}
\begin{aligned}
D(S_n|\bm{Y}_n=\bm{y}_n)\le  C(\varepsilon,s,r,\kappa)\left( \frac{M}{n^{1-2s}}+3 \right)^{1/2}\left(\frac{M^{3/k}}{n^{1/2-3s/2}}  \right)^{1/2-\varepsilon}
\end{aligned}
\end{equation}
once $n>M^{\frac{2}{kr}}.$

On the other hand, 
by condition 2 and Chebyshev inequality, we obtain
\begin{equation}\label{pro bound for A1M}
\P(A_1(M))\ge 1-\dfrac{\E J^u(X|Y)}{n^{su}},
\end{equation}
and
\begin{equation}\label{pro bound for A2M}
\P(A_2(M))\ge 1-\frac{\E|X|^k}{M}.
\end{equation}
Thus
\begin{equation}\label{pro bound for AM}
\P(A(M))\ge 1-n\left(\frac{\E|X|^k}{M}+\dfrac{\E J^u(X|Y)}{n^{su}}\right).
\end{equation}

In addition, for any $\bm{y}_n$ we use \eqref{bound for D when J large} and   \cite[Lemma 1.3]{Carlen1991} to obtain  that $\forall t>0$,
\begin{equation}\label{control D by v,J}
\begin{aligned}
D(S_n|\bm{Y}_n=\bm{y}_n)&\le C(t)(\V^t(S_n|\bm{Y}_n)+J^t(S_n|\bm{Y}_n))\\
&\le C(t)\left( \V^t(S_n|\bm{Y}_n)+\left( \frac{1}{n}\sum_{i=1}^nJ(X_i|y_i) \right)^t\right) ,
\end{aligned}
\end{equation}
where $C(t)$ is defined in \eqref{C(t)}.

Combining \eqref{bound in AM} and \eqref{control D by v,J}, we have that for $\varepsilon\in(0,1/2), s\in(1/u,1),  t\in(0,1),$ and $n>M^{\frac{2}{kr}}$,
\begin{equation}\label{fin}
\begin{aligned}
&\E D(S_n|\bm{Y}_n)=   \E D(S_n|\bm{Y}_n=\bm{y}_n)1_{\bm{y}_n\in A(M)}+\E D(S_n|\bm{Y}_n=\bm{y}_n)1_{\bm{y}_n\notin A(M)}\\
&\overset{a}{\le }  C\left( \frac{M^{4/k}}{n^{1-2s}}+3 \right)^{1/2}\left( \frac{M^{3/k}}{n^{1/2-3s/2}}  \right)^{1/2-\varepsilon}\\
&\;\;\;\;+\E C(t)\left(\V^t(S_n|\bm{Y}_n)+\left( \frac{1}{n}\sum_{i=1}^nJ(X_i|y_i) \right)^t\right)1_{\bm{y}_n\notin A(M)}\\
&\overset{b}{\le }   C\left( \frac{M^{4/k}}{n^{1-2s}}+3 \right)^{1/2}\left(  \frac{M^{3/k}}{n^{1/2-3s/2}}  \right)^{1/2-\varepsilon}\\
&\;\;\;\;+C(t)\left(\V(X)+\E J(X|Y)\right)^t\P^{1-t}\left( \bm{y}_n\notin A(M)\right)\\
&\overset{c}{\le }   C\left( \frac{M^{4/k}}{n^{1-2s}}+3 \right)^{1/2}\left(  \frac{M^{3/k}}{n^{1/2-3s/2}}  \right)^{1/2-\varepsilon}+C\left(n\left(\frac{\E|X|^k}{M}+\dfrac{\E J^u(X|Y)}{n^{su}}\right)\right)^{1-t}\\
&\overset{d}{= } O\left(\left( \frac{1}{n}\right)^{\frac{u-\frac{6}{k}u-3}{4u+\frac{6}{k}u+3}(1-2\varepsilon)}\right),
\end{aligned}
\end{equation}
where '$a$' holds due to \eqref{bound in AM} and \eqref{control D by v,J}; '$b$' comes from Cauchy inequality and $\V(X)=\E\V(X|Y)+\V\R(X|Y)$; '$c$' is deduced by  \eqref{pro bound for AM}; $d$ comes from taking
$$s=\frac{5}{3+6u/k+4u},\;\; M=n^{su},\;\; t=2\varepsilon.$$ If we take $r=su$, \eqref{fin} holds for all $n\ge1$.
Since $\varepsilon\in(0,1/2)$ is arbitrary, the proof completes.
\e
\begin{remark}
	If $X$ has the form $X=X_0+G_a$, where $G_a\sim\mathcal{N}(0,a)$ is independent of $(X_0, Y)$,  then all conditions in Theorem  \ref{rate theorem ||infty} are hold except $k$th moment condition for $X$. If we further assume the existence for any moment of $X$, then
	$$\E D(S_n|\bm{Y}_n)=O(n^{-1/4+\varepsilon}),\;\;\forall \varepsilon>0.$$
\end{remark}

\begin{remark}\label{64}
	A well-known rate theorem for central limit theorem is proved by Bobkov \textit{et al.} in \cite{Esseenbound} as follows.
	\begin{lemma}\label{Essbound}
		For independent sum $T_n=\sum_{i=1}^n \xi_i$,  $ \sum_{i=1}^n \V(\xi_i):=B_n$, if we have $\sup_i D(\xi_i)\le D<\infty$, then
		$$D(T_n)\le ce^{62D}L_4,$$
		where $c$ is a constant,
		$$D(\xi):=\KL(\xi\|\mathcal{N}(\E \xi, \V(\xi))),$$
		$$L_s:=\sum_{i=1}^n \E |\xi_i-\E \xi_i|^s/B_n^{s/2}.$$
	\end{lemma}
Compared with our Theorem \ref{KL after local limit theorem} for independent sum, Lemma \ref{Essbound} might be quiet better, if we ignore the constant influence, such as  "62" in Lemma \ref{Essbound}. However, when we try to use Lemma \ref{Essbound} to obtain convergence rate   for conditional central limit theorem, we haven't found any better way, except or assuming that conditional   variance and Fisher information satisfy that $\exists u>64$,
    \begin{equation}\label{high}
    \E\V^u(X|Y):=\int \V^u(X|Y=y)\P_Y(y)<\infty,\;\; \E J^u(X|Y)<\infty,
    \end{equation}
    to ensure the convergence (see \cite{ECCLT} for details), 
     which is quiet harsher compared with the assumption 2 in Theorem \ref{rate theorem ||infty}. 
     
     Note that the high moment assumption \eqref{high} comes from the utilization of \eqref{bound for D when J large} and truncation method. There may be a better way to utilize Lemma \ref{Essbound}, but unfortunately we haven't found it.
\end{remark}	
	
	\section*{Acknowledgments}
	\noindent We thank Zhi-Ming Ma  for his careful reading of proofs in the manuscript and helpful comments.

	\section*{Funding}
	\noindent  Yao L.-Q. was partially supported by National Key R\&D Program of China No. 2023YFA1009603.
	
		\noindent Liu S.-H. was partially supported by National Nature Science Foundation of China NSFC 12301182.

	\appendix
\setcounter{proposition}{0}
\renewcommand{\theproposition}{A.\arabic{proposition}}
\setcounter{lemma}{0}
\renewcommand{\thelemma}{A.\arabic{lemma}}

\section{Proof of $\|p_n-p_Y\|_2\to 0 \Rightarrow \|p_n-p_Y\|_1\to 0$}\label{A1}
\begin{lemma}[Exercise 4.13.1 in \cite{{functionanalysis}}]
	$\forall a,b\in\R$,
	\begin{equation}\label{a+b}
	 	||a+b|-|a|-|b||\le 2|b|.
	\end{equation}
\end{lemma}
\pf It is easy to check that $-|b|\le |a+b|-|a|\le |b|$ and \eqref{a+b} can be deduced immediately.
\e

\begin{lemma}[Exercise 4.13.2 in \cite{{functionanalysis}}]\label{Lemma A.2}
Let $(f_n)$ be a sequence in $L^2(\Omega)$ such that
\begin{enumerate}
	\item $f_n(x)\to f(x), a.e.,$
	
	\item $(f_n)$ is bounded in $L^1(\Omega)$, \textit{i.e.} $\|f_n\|_1\le M,\forall n.$
\end{enumerate}
Then $f\in L^1(\Omega)$ and 
\begin{equation}
\lim_{n\to\infty} \int (|f_n|-|f_n-f|)=\int |f|.
\end{equation}
\end{lemma}
\pf Firstly, by Fatou Lemma, $M\ge \liminf_{n\to\infty} \int |f_n|\ge \int \liminf_{n\to\infty} |f_n|=\int |f|$, thus $f\in L^1(\Omega)$. Next, let $\phi_n(x):=||f_n-|f_n-f|-|f||$, by \eqref{a+b} we have $\phi_n(x)\le |f|\in L^1(\Omega)$, thus by dominated convergence theorem,
\begin{equation}
\lim_{n\to\infty} \int \phi_n= \int \lim_{n\to\infty}\phi_n=0,
\end{equation}
and the proof completes.
\e

 \begin{lemma}[Exercise 4.13.3 in \cite{{functionanalysis}}]\label{Lemma A.3}
 	Let $(f_n)$ be a sequence in $L^2(\Omega)$ and $f$ be a function in $L^1(\Omega)$ such that
 	\begin{enumerate}
 		\item $f_n(x)\to f(x), a.e.,$
 		
 		\item $\|f_n\|_1\to \|f\|.$
 	\end{enumerate}
 	Then $\|f_n-f\|_1=0$.
 \end{lemma}
\pf Note that  $\|f_n\|_1\to \|f\|$ can deduce the assumption 2 in Lemma \ref{Lemma A.2}, thus $\lim_{n\to\infty} \int (|f_n|-|f_n-f|)=\int |f|,$  and $\lim_{n\to\infty}\int |f_n-f|=0$ is proved immediately.
\e

\textbf{Proof of $\|p_n-p_Y\|_2\to 0 \Rightarrow \|p_n-p_Y\|_1\to 0$.} We only need to prove that $\limsup_{n\to\infty} |p_n-p_Y|=0$. Let the subsequence $p_{n_k}$ satisfy $\lim_{k\to\infty} |p_{n_k}-p_Y|=\limsup_{n\to\infty} |p_n-p_Y|$, since
$p_{n_k}\overset{L^2}{\to}p_Y$, there exists a further subsequence $p_{n_{k_l}}\overset{a.e.}{\to}p_Y$ according to \cite[Theorem 4.9]{functionanalysis}. Since $\|p_{n_{k_l}}\|_1=\|p_Y\|=1$, by Lemma \ref{Lemma A.3}, we have $\lim_{l\to\infty} |p_{n_{k_l}}-p_Y|=0$, which means $\limsup_{n\to\infty} |p_n-p_Y|=0$.
\e

\section{Proof of \cite[Exercise 4.16.2]{functionanalysis} and \cite[Exercise 4.19.1]{functionanalysis}}\label{A2}

\begin{lemma}[Exercise 4.16.2 in \cite{functionanalysis}]
	Let $1<p<\infty$ and $(f_n)$ be a sequence in $L^p$ such that
	\begin{enumerate}
		\item $f_n$ is bounded in $L^p$.
		
		\item $\|f_n-f\|_1\to 0$.
	\end{enumerate}
Then $f_n\rightharpoonup f$ weakly $\sigma(L^p, L^{p'})$.
\end{lemma}
\pf On the one hand,  since $L^p$ is reflexive for any p, $1<p<\infty$, for any subsequence of $f_n$,  there exists a further subsequence $f_{n_k}$ and $g\in L^p$ such that $f_{n_k}\rightharpoonup g$ weakly $\sigma(L^p, L^{p'})$. On the other hand, since $\|f_{n_k}-f\|_1\to 0$, there exists a further subsequence $f_{n_{k_l}}$ such that $f_{n_{k_l}}\xrightarrow{a.e.} f$, therefore $f=g$ (by the \textbf{Hint}). In conclusion, any subsequence of $f_n$ has a further subsequence  weakly convergence to $f$, thus  $f_n\rightharpoonup f$ weakly $\sigma(L^p, L^{p'})$. 
\e

\begin{lemma}[Exercise 4.19.1 in \cite{functionanalysis}]
	Let $1<p<\infty$, $(f_n)$ be a sequence in $L^p$ and $f\in L^p$ such that
	\begin{enumerate}
		\item $f_n\rightharpoonup f$ weakly $\sigma(L^p, L^{p'})$.
		
		\item $\|f_n\|_p\to\|f\|_p$.
	\end{enumerate}
	Then $f_n\xrightarrow{L^p}f$.
\end{lemma}
\pf Argue by contradiction. 
If there exists a subsequence that $\lim_{k\to\infty}\|f_{n_k}-f\|_p=a>0$, since $f_{n_k}$ is bounded, there exists a further subsequence $f_{n_{k_l}}$ and $g\in L^2$ such that $f_{n_{k_l}}\xrightarrow{L^p} g$ and $f_{n_{k_l}}\rightharpoonup  g$ weakly $\sigma(L^p, L^{p'})$. Therefore $f=g$ and $\lim_{l\to\infty}\|f_{n_{k_l}}-f\|_p=0$. A contradiction!
\e


\end{document}